# Greedy Placement of Measurement Devices on Distribution Grids based on Enhanced Distflow State Estimation

Paola Paruta, Thomas Pidancier, Mokhtar Bozorg*, Mauro Carpita

*School of Management and Engineering Vaud (HEIG-VD)*
*University of Applied Sciences of Western Switzerland (HES-SO)*
Route Cheseaux 1, Case Postale 521, 1401 Yverdons-les-Bains, Switzerland
{paola.paruta,Thomas.pidancier,mokhtar.bozorg,mauro.carpita}@heig-vd.ch

* Corresponding author: mokhtar.bozorg@heig-vd.ch

*Abstract*— The needs for improving observability of medium and low voltage distribution networks has been significantly increased, in recent year. In this paper, we focus on practical approaches for placement of affordable Measurement Devices (MDs), which are providing three phases voltage, current, and power measurements with certain level of precision. The placement procedure is composed of a state-estimation algorithm and of a greedy placement scheme. The proposed state-estimation algorithm is based on the Distflow model, enhanced to consider the shunt elements (e.g., cable capacitances) of the network, which are not negligible in low voltage networks with underground cables. The greedy placement scheme is formulated such that it finds the location of minimum required number of MDs while certain grid observability limits are satisfied. These limits are defined as the accuracy of state-estimation results in terms of voltage magnitudes and line currents over all nodes and lines, respectively. The effectiveness of the proposed placement procedure has been validated on a realistic test grid of 10 medium voltage nodes and 75 low voltage nodes, whose topology and parameters were made available from the Distribution System Operator (DSO) of the city of Geneva, Switzerland.

*Index Terms*— Distflow, Distribution network monitoring, greedy algorithms, Measurement devices placement, State estimation.

## I. INTRODUCTION

IN recent years, the uncontrolled injection of renewable energy (e.g. PV rooftop installations) on the electrical distribution network, as well as the electrification of the automotive sector, has contributed to the destabilization of medium and low voltage grids. Thus, the needs for network monitoring and control solutions has been significantly increased. Regulators, all around the world, encourages the Distribution System Operators (DSO) to improve the network observability by installing low voltage grid Measurement Devices (MDs) as well as smart-meters (SM) on the client sides. For instance, in Switzerland, at least 80% of the DSO's clients should be equipped with smart-meters by the end of 2027 [1]. Nevertheless, medium and low voltage distribution networks are very large, in terms of number of nodes and lines. Hence, the DSOs are looking for affordable and scalable solutions for the network supervision, providing them easy access to the network statistics, power quality and security status.

In this paper, we focus on practical approaches for placement of affordable measurement devices, which are providing three phases voltage, current, and power measurements with certain level of precision. The placement procedure developed to determine the required number and location of measurement devices, is composed of a state-estimation algorithm and of a greedy placement scheme.

Recently, a number of Distribution System State Estimation (DSSE) methods are proposed in the literature. The methods differ in terms of; i) available measurement data (e.g., availability of synchronous measurements including magnitude and angle of voltages and currents [2]), ii) state variables (branch-current, and node-voltage), iii) network topology and line models, and iv) state estimation algorithm (e.g., Weighted Least Square estimators and iterated Kalman-filter [3]). A comprehensive review of the state-of-the-art methods and techniques is presented in [4]. State estimation algorithms in most of the presented methods, are developed based on both voltage angles and voltage magnitudes. However, voltage angles are not always available with affordable measurement devices. Note that in a Low Voltage (LV) distribution network, the lines are relatively short. Therefore, the difference between voltage angles of neighboring nodes are very low and can be in the range of the measurement precision.

In [5], a state estimation method based on Distflow model [6] is introduced, which does not require measurement of voltage angles, but the shunt element (capacitance) of the lines are not modeled.

In this paper, the proposed state-estimation algorithm is based on the Distflow model, enhanced to consider the shunt elements of the network using a PI model of the lines. Note that the cables' capacitance is not negligible in low voltage networks with underground cables. Indeed, the enhanced method is well suited for radial low voltage distribution networks where the measurement of voltage angles with high precision is not always available. It uses the Weighted Least Square method with a simplified constant Jacobian matrix derived from the Distflow equations.

The algorithm estimates the state of the grid by using; i) voltage magnitude and power flow measurements (10-minutes average values) coming from MDs installed at specified nodes, and ii) load pseudo-measurements (PMs) extrapolated from day-ahead forecast using aggregated SM data (active and reactive power consumption of end-consumers connected to a given node), at every node. This is actually a day-ahead forecast, as the SM data are not available in real-time, due to legal and privacy constraints.

The placement scheme is formulated such that it finds the location of minimum required number of MDs while certain grid observability limits are satisfied. These limits are defined as the accuracy of state-estimation results in terms of voltage magnitudes and line currents over all nodes and lines, respectively. The algorithm we use to optimize the placement of MDs follows a greedy approach. Greedy algorithms are a class of algorithms that aim at finding the global solution of an optimization problem by choosing the local optimal solution at each iteration. Although this approach does not guarantee to find the global optimum solution, it has been widely used in power system applications (e.g., Placement of phasor measurement units in transmissions networks [7], observability analysis [8], and reconfiguration of medium voltage distribution networks [9]) thanks to its tractability and applicability.

In the proposed greedy placement scheme, at each iteration, we add one measurement device. The location of the new measurement device is found such that a state-estimation accuracy function in presence of existing MDs (from previous iterations) and the new MD is maximized. The iterative process will be terminated as soon as grid observability limits are satisfied.

Finally, we applied the proposed placement procedure on an 85-nodes simulation case derived from a medium and low voltage distribution grid in Geneva, Switzerland. The numerical results demonstrate that the proposed approach is tractable and computationally efficient.

The rest of the paper is organized as follows. The enhanced Distflow state-estimation method is presented in section II. The iterative greedy placement scheme is formulated and discussed in section III. Numerical results and assessment of the proposed methods are discussed in section IV. Finally, conclusions are presented in section V.

## II. DISTRIBUTION SYSTEM STATE-ESTIMATION

### A. Hypothesis and context

In order to define and formulate the proposed state estimation algorithm, the following is assumed:

- The distribution network configuration is radial, and modelled by load-buses and a slack-bus. The secondary side of the high-voltage (HV) / medium-voltage (MV) transformer is considered as the slack-bus.
- Load pseudo-measures (active and reactive power) are available at every nodes of the network. Their values can be negative (injection), null (connection node) or positive (consumption).
- The MD installed at any node, gives the voltage magnitude and the active and reactive power flow in all the lines/transformers connected to the node.
- The lines and transformers parameters (i.e., rated values and impedances) are known.

Figure 1 illustrates the hypotheses on the topology and the available measurements.

Assuming pseudo-measurement (PM) of the power load is available at every node and under the assumption that the line parameters are known, only one voltage measurement is needed to fully determine the state of the grid. Indeed, since every PM provides information of the power load at its node, once all of the power loads are known, the state of the grid can be simply computed, using the Distflow algorithm. The need to install additional MDs comes from the fact that PMs might not be sufficiently accurate, as the SM data are not available in real-time due to legal and privacy issues. For instance, in Switzerland, the SM data will only be available to the DSO every 24h, with a granularity of 15min, for privacy reasons. Between data deliveries, some forecasting will be needed to mimic real-time signals. Additionally, the SM data needs to be aggregated, as the forecast power consumption per node can be derived from a number of buildings (SMs) connected to that node. The PMs will be necessarily affected by an error, which will impair the accuracy of the predicted state. The role of the MDs is to compensate such uncertainties by providing real-time, precise measurements of power flows, voltages and currents through the nodes/lines where they are installed.

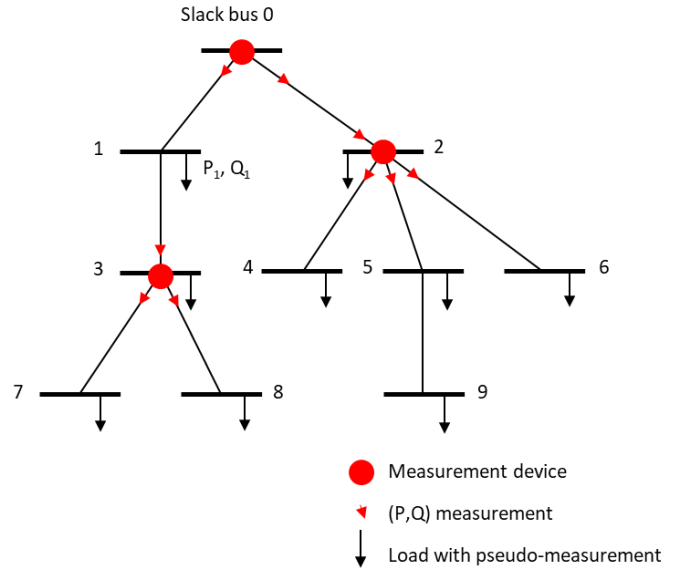

Figure 1. Network example

### B. Distflow algorithm

The Distflow algorithm [6] is a load-flow method that follows an iterative pattern well-suited for radial networks and does not rely on the voltage angle. It can also be generalized to any network topology [10] and work with PI model of lines (Figure 2). As discussed in II.A, the studied network is assumed radial, with a slack-bus at the transformation node, $L$ lines and $L$ load-buses nodes.

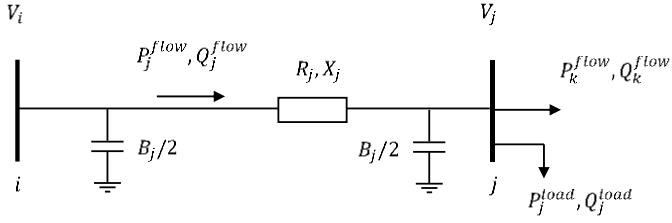

Figure 2. Pi model of a line and notations used in the Distflow

The variables are computed in two steps, first the power flow in the lines from the bottom-up is computed using the backward-distflow Equation (1), then the square voltage magnitude of the nodes from the slack bus to the end-nodes is derived using the forward-distflow Equation (2). The capacitance of the lines is considered in the computation of the reactive power flow. These calculations are repeated until the convergence is obtained.

$$\begin{cases} P_j^{flow} = P'_j + R_j \frac{P'^2_j + Q'^2_j}{V_j^2} \\ Q_j^{flow} = Q'_j + X_j \frac{P'^2_j + Q'^2_j}{V_j^2} \end{cases} \quad (1)$$

with $\begin{cases} P'_j = P_j^{load} + \sum_k P_k^{flow} \\ Q'_j = Q_j^{load} + \sum_k Q_k^{flow} + Q_{j_B} \end{cases}$ where $Q_{j_B} = -B_{\Sigma j} V_j^2$

As indicated in Figure 2, $P_j^{flow}$ and $Q_j^{flow}$ are active/reactive power flow in line $j$ (which is upstream to node $j$). $P_j^{load}$ and $Q_j^{load}$ are active/reactive power load at node $j$. $V_j$ is the square of voltage magnitude at node $j$. $R_j$ and $X_j$, are resistance and reactance of line $j$, respectively. The term $B_{\Sigma j}$ presents sum of the shunt elements (susceptances) of all the lines connected to node $j$.

$$V_j^2 = V_i^2 - 2(R_j P_j^{flow} + X_j Q_j^{flow}) + (R_j^2 + X_j^2) \frac{P_j^{flow\,2} + Q_j^{flow\,2}}{V_i^2} \quad (2)$$

### C. Enhanced distflow state estimation

The developed state estimation algorithm is based on the Distflow model presented in section II.B, that takes into account the susceptances (shunt capacitance) of the lines. This is an improvement of the original method proposed in [5]. It uses the weighted-least-square (WLS) estimation method and a set of measurements $z$ (3).

$$z = h(x) + e \quad (3)$$

Where $x$ is the state vector, $h$ the measurement function between $x$ and $z$, and $e$ is the measurement noise. Following the definition of the enhanced Distflow algorithm, we see in (1) and (2) that it is possible to define the state estimation vector by (4).

$$x = \begin{bmatrix} P^{load} \\ Q^{load} \\ V_0^2 \end{bmatrix} \quad (4)$$

$P^{load} = [P_1^{load} \dots P_L^{load}]^T$ and $Q^{load} = [Q_1^{load} \dots Q_L^{load}]^T$ are the active and reactive power load of every load-buses, respectively. $V_0^2$ is the square of voltage magnitude at slack-bus.

The pseudo-measurements (derived from day ahead forecasts) are available at every node ($\widehat{P}^{load}, \widehat{Q}^{load}$), and the MDs placed in the network measure the node voltage and the power flow in the connected lines ($\widehat{V}, \widehat{P}^{flow}, \widehat{Q}^{flow}$). Hence, the measurement set $z$, can be written by (5).

$$z = \begin{bmatrix} \widehat{P}^{load} \\ \widehat{Q}^{load} \\ \widehat{P}^{flow} \\ \widehat{Q}^{flow} \\ \widehat{V}^2 \end{bmatrix} \quad (5)$$

It is possible to solve the state estimation problem in an iterative manner with the help of (6),

$$x_{(k+1)} = x_{(k)} + (H^t W H)^{-1} H^t W r_{(k)} \quad (6)$$

where $k$ is the iteration index. $H$ is the jacobian matrix of $h$ as formulated in (7).

$$H(x) = \begin{bmatrix} \frac{\partial P^{load}}{\partial P^{load}} & \frac{\partial P^{load}}{\partial Q^{load}} & \frac{\partial P^{load}}{\partial V_0^2} \\ \frac{\partial Q^{load}}{\partial P^{load}} & \frac{\partial Q^{load}}{\partial Q^{load}} & \frac{\partial Q^{load}}{\partial V_0^2} \\ \frac{\partial P^{flow}}{\partial P^{load}} & \frac{\partial P^{flow}}{\partial Q^{load}} & \frac{\partial P^{flow}}{\partial V_0^2} \\ \frac{\partial Q^{flow}}{\partial P^{load}} & \frac{\partial Q^{flow}}{\partial Q^{load}} & \frac{\partial Q^{flow}}{\partial V_0^2} \\ \frac{\partial V^2}{\partial P^{load}} & \frac{\partial V^2}{\partial Q^{load}} & \frac{\partial V^2}{\partial V_0^2} \end{bmatrix} \quad (7)$$

$$W = \begin{bmatrix} \sigma_{z_1}^2 & \cdots & 0 \\ \vdots & \ddots & \vdots \\ 0 & \cdots & \sigma_{z_m}^2 \end{bmatrix}^{-1} \quad (8)$$

$W$ is the inverse of the variance matrix formulated in (8) regarding all the measurements, where $\sigma_{z_n}$ is the standard deviation of the $n^{th}$ measurement. The procedure for computation of these measurement uncertainties ($\sigma_{z_n}$) is discussed in the following section III.A.

$r_{(k)} = z - h(x_{(k)})$ is the error between the measurement and the results of the application of the measurement function to the computed state vector at iteration $k$.

We see in (7) that $H$ depends on the state of the system, but with a series of approximations as introduced in [5], it is possible to give a fix value to $H$ that depends only on the network's characteristics (e.g., line resistance, reactance, etc.). The different elements of $H$ are set as shown in (9)-(14). Equation (12) shows the relation between the capacitance of the lines and the reactive power flow.

$$\frac{\partial P^{load}}{\partial P^{load}} = \frac{\partial Q^{load}}{\partial Q^{load}} = 1 \quad (9)$$

$$\frac{\partial P^{load}}{\partial V_0^2} = \frac{\partial Q^{load}}{\partial V_0^2} = \frac{\partial P^{flow}}{\partial V_0^2} = \frac{\partial P^{flow}}{\partial Q^{load}} = \frac{\partial Q^{flow}}{\partial P^{load}} = 0 \quad (10)$$

$$\frac{\partial P^{flow}_i}{\partial P^{load}_j} = \frac{\partial Q^{flow}_i}{\partial Q^{load}_j} = \begin{vmatrix} 0 \\ or \\ 1 \end{vmatrix} \quad (11)$$

1 if the node $j$ is after the line $i$, 0 otherwise

$$\frac{\partial Q^{flow}_i}{\partial v_0^2} = -\sum_j \frac{B_j}{2} \quad (12)$$

with $j = \{index\ for\ all\ the\ lines\ after\ i\}$

$$\frac{\partial V_i^2}{\partial P^{load}_j} = -2\sum_k R_k \quad (13)$$

$$\frac{\partial V_i^2}{\partial Q^{load}_j} = -2\sum_k X_k \quad (14)$$

where $k$ is the index of the lines between the node 0 and the furthest node that leads to $i$ and $j$, as indicated in see Figure 3.

$$\frac{\partial V_i^2}{\partial V_0^2} = 1 \quad (15)$$

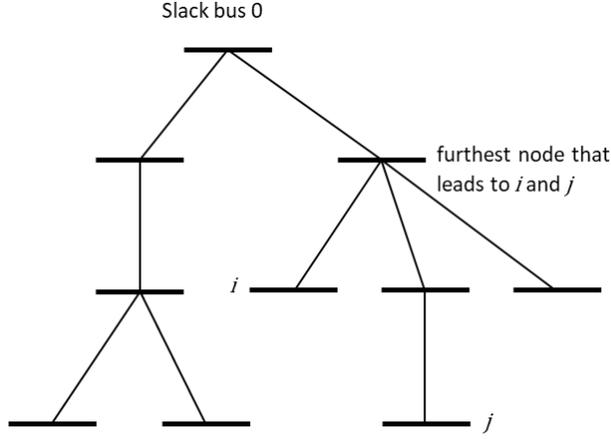

Figure 3. Example on how the coefficient $\frac{\partial V_i^2}{\partial P^{load}_j}$ and $\frac{\partial V_i^2}{\partial Q^{load}_j}$ are calculated

Equations (9)-(14) approximates the component of the Jacobian matrix $H$. With this constant approximated $H$, equation (6) can be solved iteratively and eventually converge to solution $x^s$. Afterwards, the estimated state of the grid $\chi$ as defined in (16), is computed by running the iterative Distflow algorithm (1)-(2), with $x^s$ as input.

$$\chi = \{\{P^{load}, Q^{load}, V\}_i, \{P^{flow}, Q^{flow}\}_l\} \quad (16)$$

Where $i = [\![1, L]\!]$, and $l = [\![1, L]\!]$ represents index of nodes and lines of the network, respectively.

### III. GREEDY PLACEMENT SCHEME

The greedy placement scheme presented hereafter shows a way of scanning through all the possible sets of measurement devices (MDs) locations, to find the one that allows to estimate the grid state with the desired accuracy, while aiming to measure the least number of nodes.

In this section, first, the uncertainties associated with pseudo-measurements (mainly due to day-ahead forecast error) as well as measurement devices are formulated. Then, a cost function that represent the inaccuracy of the results of the state estimation algorithm using a predefined set of measurement devices is introduced. Afterwards, the iterative greedy approach is presented in details.

#### A. Measurement and pseudo-measurements uncertainties

Let's define a grid state $\chi^*$ as follows:

$$\chi^* = \{\{P^{load*}, Q^{load*}, V^*\}_i, \{P^{flow*}, Q^{flow*}\}_l\} \quad (17)$$

Concretely, in a real distribution network, the grid state $\chi^*$ is unknown. Here we introduce it only for simulation purposes, in order to generate a coherent set of PMs and MDs.

We emulate the uncertainty of the pseudo-measurements at each node $i$, by extracting the value of the PM from a normal distribution centered around $\chi_i^{PM*} = \{P^{load*}, Q^{load*}\}_i$:

$$z_i^{PM} \sim N(\chi_i^{PM*}, \sigma_i^{PM}) \quad (18)$$

The standard deviation $\sigma_i^{PM}$ will depend on the uncertainty level of the PM at node $i$. In general, the higher the power load, the higher the uncertainty of the PM:

$$\sigma_i^{PM} = C^{PM}\chi_i^{PM*} + \sigma_0^{PM} \quad (19)$$

where $\sigma_0^{PM}$ is a small absolute error and $C^{PM}$ represents the constant of proportionality. Regarding the error of day-ahead nodal forecasting algorithms, in the numerical simulations, $C^{SM}$ is set as 0.2 (20% uncertainty) and $\sigma_0$ is set as $10^{-4}$ (0.01%).

In this work, the selection of normal distribution is justified as day-ahead forecast of load at each node considers aggregated data from a number of buildings (smart meters) connected to that node with different load patterns. However, one could argue that certain load components such as electrical heating are correlated, and underestimating the temperature in the forecasting model would lead to a correlated forecast error across the loads of the network.

Measurement device precision level is used to model uncertainties associated with measurement data, as:

$$z_{i,l}^{MD} \sim N(\chi_{i,l}^{MD*}, \sigma_{i,l}^{MD}) \quad (20)$$

where $\chi_{i,l}^{MD*} = \{\{V^*\}_i, \{P^{flow*}, Q^{flow*}\}_l\}$ for $(i, l)$ being the index of nodes/lines where MDs are installed. Like pseudo-measurement, the standard deviation is a linear function of the measurement as

$$\sigma_{i,l}^{MD} = C^{MD}\chi_{i,l}^{MD*} + \sigma_0^{MD} \quad (21)$$

Note that the uncertainty for the MD data is much smaller than the uncertainty on a PM forecast, i.e. $C^{MD} \ll C^{PM}$. In the numerical simulations we set $C^{MD} = 0.005$ (0.5 %) as uncertainty on $P^{flow}, Q^{flow}$ and $C^{MD} = 0.001$ (0.1 %) for measurements of voltage magnitude square $V^2$.

#### B. Inaccuracy cost function: evaluation of a placement configuration

In this section, we present how to compute the inaccuracy cost associated to a predefined set of installed measurement devices. We refer to this set as the MDs configuration with $g$. For example, if measurements are placed at node 0, 5, and 14 on a grid, then $g = \{0, 5, 14\}$. Given a grid state $\chi^*$, in order to determine the cost $J(g)$, the algorithm randomly generates, through Equations (18, 20), a vector of MDs and PMs measurements $z_g = \{z_g^{MD}, z^{PM}\}$ and solves the state estimation to obtain $\chi$. Depending on how many measurements we placed and their accuracy, $\bar{\chi}$ will be more, or less, close to the initial state. Hence, the difference between the estimated state $\chi$ and the actual state $\chi^*$ is a first indicator of how good the choice of $g$ is. Nonetheless, this difference is also impacted by the initial random guess through which the measurement vector was generated. To alleviate such dependency, the process is repeated for a sufficiently large number of times. For each new generated vector $z_g^{[r]}$, the state $\chi^{[r]}$ is determined, with the apex $[r]$ ($r = 1, ..., R$) indicating a random realization. In the following, we select $R = 20'000$, so that the relative difference between the standard deviation computed from the random realizations matches, and the one chosen as input is less than 2%. The average of all the $\chi^{[r]}$ is the exact state $\chi^*$, but the variance of the distribution of all $\chi^{[r]}$ is dependent on $g$, where a smaller variance indicates a better choice of location $g$. Through the same reasoning detailed in section III.A, this variance can be reinterpreted as the uncertainty that one can expect on the estimation of the grid state, given a set of locations $g$, $\sigma_g$ such that $\chi \sim N(\chi^*, \sigma_g)$.

The smaller the variance, the higher the accuracy.

Finally, the cost vector $J(g)$ we use to frame the optimization problem is defined as:

$$J_i(g) = \max(\sigma_{g,i} - \sigma_{max,i}, 0) \quad (22)$$

with $\sigma_{max,i}$ the maximal accepted uncertainty that the DSO requires at node/line $i$.

Equation (22) means that the cost vector $J$ of a set of placements $g$ is the uncertainty that exceeds the desired maximal one. If the uncertainty $\sigma_{g,i}$ at a node/line $i$ is lower than the desired one, than that node/line does not bring any additional cost. Note that if $||J(g)||_\infty = 0$ than the grid state is predicted with the desired accuracy.

In practical applications the desired accuracy $\sigma_{max,i}$ can be formulated so that the threshold is set on the relative uncertainty rather then the absolute one, for example $\sigma_{max,l} = 0.1 P_l^{cap}$, if we want the uncertainty on the power flow on a line, normalized to the line capacity $P^{cap}$, to be below 10%, i.e. $\sigma_g(P_l^{flow})/P^{cap} < 0.1$.

*C. Objective function*

In section III.B we defined the inaccuracy cost function, with $g$, the set of MDs placements as decision variable. The objective of the optimization algorithm is to find the minimum number of MDs that allow observability of the grid with the desired accuracy:

$$\text{Objective:} \min_{g:||J(g)||_\infty=0} |g| \quad (23)$$

where $|g|$ indicates the number of MDs.

*D. Finding the optimal placement*

Inspired by the iterative placement method proposed in [11], the algorithm we developed to optimize the placement of MDs follows a greedy iterative approach. Each iteration $t$ results in the addition of a measurement with respect to the configuration at $t-1$, so that $|g|^t = |g|^{t-1}+1$. The location of the new measurement, $b^t$, is chosen between the ones available in order to minimize the cost function:

$$b^t = \min_{i \in A^t} ||J(g^{t-1} \cup \{i\})||_\infty \quad (24)$$

where $J$ is the cost vector function in Equation (23), and $A^t$ is the set of available nodes at iteration $t$. $A^t$ contains all nodes that do not already have a MD installed. The *leaves* are also excluded from $A^t$. Following the tree analogy, these are all the nodes that only have an incoming line. We discard them because placing a MD on a *leaf* is never better than placing it on the node just above. Let us recall that every MD measures $P^{flow}$ on incoming and outgoing lines as well as the voltage at the node where it is installed. Therefore, placing a MD on a node just above a *leaf* gives the $P^{flow}$ value on the line connecting to the *leaf*, the voltage at the node, and, consequently, allows the algorithm to determine exactly the voltage value at the *leaf*, too (the line properties are assumed to be known exactly). Besides, it also provides $P^{flow}$ on all connected lines. On the other hand, placing the MD on the *leaf* would have given only $V_{leaf}$ and the $P^{flow}$ on the one incoming line. In the simulation results (section IV) we show that the removal of all *leaves* as candidates for placement decreased the computational cost by a factor two.

The algorithm stops when $||J(g)||_\infty = 0$, meaning that the desired accuracy is guaranteed everywhere on the grid.

The underlying hypothesis to assure that the solution obtained is the global optimum is that the problem has "optimal substructure", i.e. that the optimum solution is found by selecting the optimum at every sub-step. Although the presented workflow is intuitive and has already been employed for the optimal placement of MDs in [11], it is not possible to prove that the solution obtained with the greedy algorithm satisfies the objective as stated in Equation (23), nonetheless, the chosen set of MDs does satisfy $||J(g)||_\infty = 0$.

*E. Grid loading scenario*

As discussed earlier, a change in MDs precision or PM accuracy can impact the placement of MDs. What is maybe less obvious is that the grid state $\chi^*$, and more precisely the grid loading scenario $x^*$, also affects the choice of optimal placement. In real application, the DSOs are interested to create worst-case scenarios where the grid is pushed to its limit, and see which nodes in that case need to be closely monitored. Note that the proposed greedy placement algorithm can take as input any arbitrary combination of loads, $P^{load}$ and $Q^{load}$, as long as the capacity (i.e., maximum current limit) of the lines is not exceeded.

Finally, it is worth mentioning that in this paper, the placement scheme is formulated for a balanced system for the sake of clarity of presentation. In reality, loading of low voltage distribution networks are unbalanced. To deal with unbalanced loading scenarios, the proposed state estimation algorithm formulated in section II, should be solved for each phase separately. The underlying assumption is that the line parameters are balanced. Afterward, in the greedy placement algorithm, the choice of state estimation uncertainty index could be defined as the worst error among three phases.

IV. SIMULATIONS AND NUMERICAL RESULTS

*A. Greedy placement on Geneva subgrid*

In this section, the proposed greedy placement method is applied on an 85-nodes simulation case derived from a medium and low voltage distribution grid in Geneva, Switzerland. The network composed of 10 medium voltage (i.e., $N_{MV} = \{0,1,2,4,10,14,19,25,28,31\}$) and 75 low voltage nodes. The grid loading scenario $x^*$, and the corresponding grid state $\chi^*$, obtained after solving the Distflow, are represented in Figure 4. The color-scheme for the nodes is based on the power load consumption $S$ in p.u. (the reference power being 10 MVA). Given the power loads everywhere, and the voltage value at node 0, the Distflow algorithm is used to compute power flows and currents in the lines. The color-scheme for the lines in Figure 4 shows the line loading $I^{flow}/I^{cap}$ in percentage, where

$$I^{flow} = \sqrt{P^{flow^2} + Q^{flow^2}}/V_i$$

and $I^{cap}$ is the maximum current limit. In this grid loading scenario, lines ending at node 84, 32, 43, 54, 65 and 52 have the highest line loading (close to 80%). Whereas, the medium voltage cables (lines connecting MV nodes), are weakly exploited between 10% to 35%.

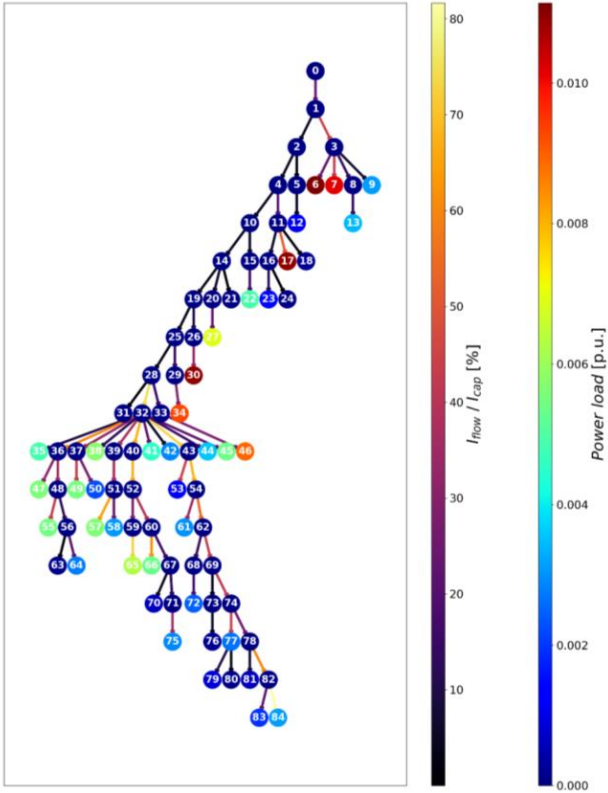

Figure 4 Grid state chosen for the 85-node sub-grid of Geneva. The color-scheme of the lines (left-scale) shows the current flow relative to the line capacity, while the nodes color-scheme (right-scale) represent the power load in p.u

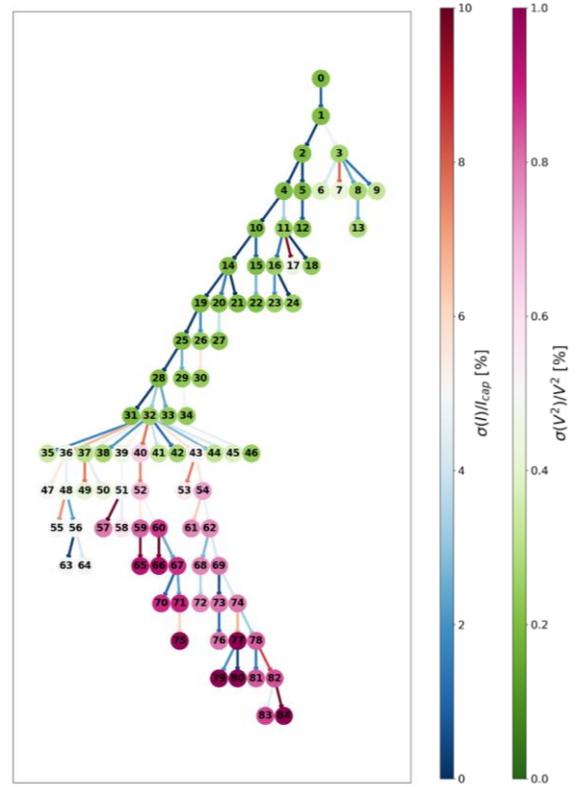

Figure 5. Results of the enhanced state estimation algorithm on the 85-node subgrid of Geneva, with grid state of Figure 4, no MD installed, only SMs at all nodes and $V_0$ at slack-bus are available

Let us define the base case, where no MD is installed, and only PM forecasted data are used, with uncertainty $C = 0.2$ and $\sigma_0 = 1.9 \times 10^{-7}$ in Equation (18) together with $V_0$. The resulting uncertainty on the grid state estimation is plotted in Figure 5.

The color-scheme shows the uncertainty the DSO can expected when solving the power-flow with the PM data and a voltage measurement at the slack bus. The nodes color-scheme is green to fuchsia for a relative uncertainty on the voltages, $\sigma(V^2)/V^2$ ranging from 0 to 1%. Nodes higher up in the tree have better accuracy, because they are closer to the measured slack-bus. The line color-scheme, blue to red, maps the relative uncertainty on current flow along the line, relative to the line capacity, $\sigma(I^{flow})/I^{cap}$ from 0 to 10%.

The maximum state-estimation uncertainties allowed by the DSO are defined as 0.3% and 5% for estimation of voltage and current square magnitudes, respectively. Therefore, with reference to Equation (22), we have $\sigma_{max,i} = 0.003 V_i^2$, and $\sigma_{max,l} = 0.05 I_l^{cap^2}$. As we can see in Figure 5 a), these limits are not satisfied in 51 nodal voltages and 23 line currents.

Next, we run the proposed greedy placement scheme to find the optimal placement of MDs, in order to satisfy the maximum state-estimation uncertainty limits. We assumed the MDs to have a relative precision of 0.2% for $V^2$ ($C = 0.002$ and $\sigma_0 = 0$ in Equation (20)) and of 0.5% for $P^{flow}$ and $Q^{flow}$ ( $= 0.005$ and $\sigma_0 = 1.9 \times 10^{-7}$). These precision levels correspond with affordable measurement devices which are commercially available.

Figures 6 a) and b) show the result of the greedy placement algorithm when optimizing either the voltage accuracy (panel a), or current accuracy (panel b). The nodes circled in orange are those where the MDs must be placed.

In this grid-loading scenario, installation of five MDs is sufficient to achieve the desired accuracy on $V^2$. Whereas, we need to install eight MDs to achieve the desired accuracy level for line currents. Node 32, which has the highest "connectivity" because 11 lines directly depart from it, is chosen by the algorithm to be measured in both cases. The placement configuration in Figure 6 a), although obtained for achieving the voltage accuracy, it already greatly improves the current accuracy. However, on a few lines the uncertainty on the estimation of current flow still exceeds $\sigma_{max} = 5\% \, I^{cap}$. In Figure 6 b), we can see that eight measurement devices are required to ensure that the uncertainty of current estimation at every line is below the allowable range $\sigma_{max} = 5\% \, I^{cap}$.

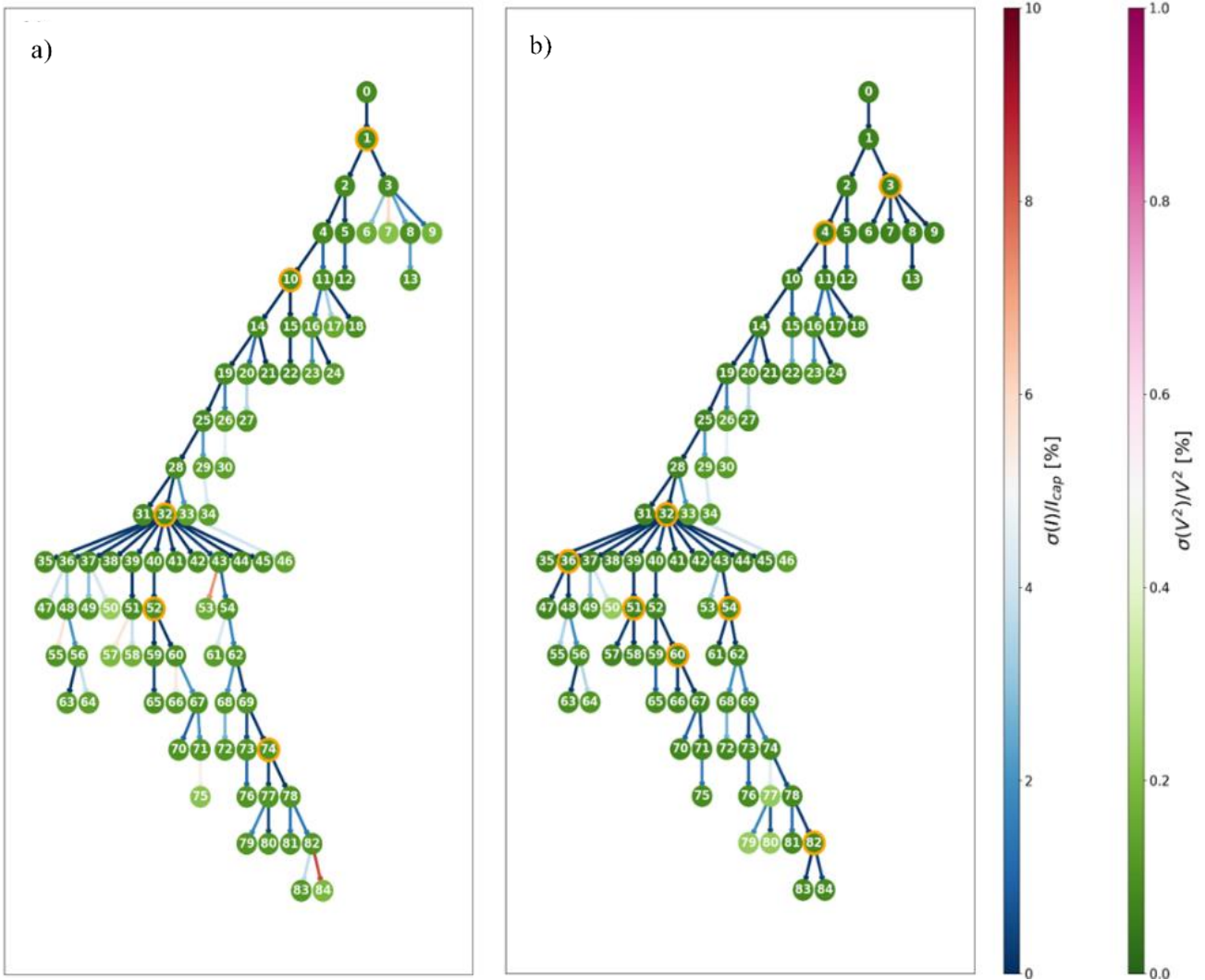

Figure 6. Greedy algorithm result for the placement of measurement devices on the 85-node subgrid of Geneva, with grid state of Figure 4, a) uncertainty after optimizing based on voltage, $\sigma(V^2)/V^2 < 0.3\%$, b) uncertainty after optimizing the placement based on current, $\sigma(I^{flow})/I^{cap} < 5\%$.

\* The orange circle around the node indicates the placement of a Measurement Devices (MD).

*B. Computational cost*

The workflow of the greedy placement algorithm described in Section D implies a computational cost and a run-time $T$ that scales as $|g| \times n$:

$$T = K[n + (n-1) + \cdots (n - (|g|-1))] = \mathcal{O}(|g|n) \quad (25)$$

where $n$ is the number of nodes available for the placement of a measurement, $|g|$ is the final number of measurements placed, and $K$ is the time needed to evaluate the cost vector $\mathbf{J}$ associated to one placement configuration, requiring the resolution of $10^3$ state estimation problems (see Section B). Equation (25) is explained by considering that at every iteration $t$, all available nodes are tested by running the state estimation with MDs placed at $g^{t-1} \cup \{i\}$. The node leading to the lower cost is then chosen and removed from the available nodes.

Concretely, for the optimal placement based on voltage accuracy in Figure 5b), the algorithm had to test $|g|n = 5 \times 42$ possible configurations. Here, $n$ is 42 rather than 85 because all 43 *leaves* have been discarded. Running the algorithm in MATLAB R2019a on a single processor (Intel® Core™ i7-4770 CPU @ 3.40GHz), leads to K = 15sec and $T \sim 40$ min. In future and for application on larger grids, the algorithm could be easily parallelized on $n$ processor and the run-time reduced by a factor $n$ to a few minutes ($T \sim K|g|$).

*C. Sensitvity analysis: Number of measurement devices needed*

We expect the number of measurements $|g|$ decreases as the allowed uncertainty $\sigma_{max}$ increases. Figure 7 shows that the algorithm for the optimal placement reproduces this expected qualitative behavior. Nonetheless, from a quantitative point of view, the exact value of $|g|$ is affected by a number of other parameters, such as the grid loading scenario or the PMs uncertainty.

Figure 7 shows that a decrease in the uncertainty of the PMs of 5%, from $C^{PM} = 0.2$ (blue line) to $C^{PM} = 0.15$ (orange line),

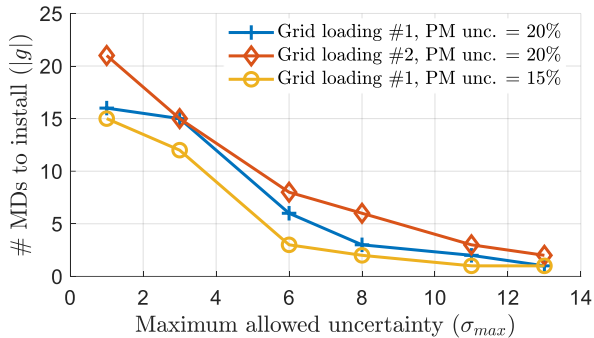

Figure 7. Number of measurement devices to install decreases when we allow for a higher uncertainty on the prediction of the grid state. The different lines correspond to different grid loading scenario and pseudo-measurement ucertainty.

then a lower number of MDs is needed. Reduction of the uncertainty of the PMs can be achieved by either improving the performance of day-ahead forecast methods for aggregated smart-meter data or making these data available near real-time.

Similarly, two different grid loading can require a different number of MDs for the same desired accuracy. In Figure 7 the grid loading scenario #1 (blue line) demands 16 MDs for an uncertainty of 1% while scenario #2 requires 21 of them (red line). An intuitive understanding of the relevance of the grid loading scenario was already discussed in Section III.E.

## V. SUMMARY AND CONCLUSIONS

In this paper, we developed a practical procedure for placement of affordable Measurement Devices (MDs), which are providing three phases voltage, current, and power measurements with certain level of precision. The placement procedure is composed of a state-estimation algorithm and of a greedy placement scheme.

The proposed state-estimation algorithm is based on the Distflow model that does not require voltage angle measurements. The model is enhanced to consider the shunt elements (e.g., cable capacitances) of the network, which are not negligible in medium and low voltage distribution networks with underground cables.

The greedy placement scheme is formulated such that it finds the location of minimum required number of MDs while certain grid observability limits are satisfied. These limits are defined as the accuracy of state-estimation results in terms of voltage magnitudes and line currents over all nodes and lines, respectively. The underlying hypothesis to assure that the solution obtained is the global optimum is that the problem has "optimal substructure", i.e. that the optimum solution is found by selecting the optimum at every sub-step. Although the presented workflow is intuitive, it is not possible to prove that the solution obtained with the greedy algorithm satisfies the global optimality condition. Nonetheless, the obtained solution (chosen set of MDs) does satisfy the desired grid observability condition (i.e., desired accuracy of state estimation).

The effectiveness of the proposed algorithm has been validated on a realistic 85-nodes distribution grid, in terms of optimality, and tractability under different desirable observability. Finally, sensitivity analysis has shown the importance of grid loading scenarios as well as the uncertainties associated with pseudo-measurements on the optimal number and locations of measurement devices. Finally, it is worth mentioning that further works are needed to evaluate the influence of the hypothesis "the line properties are assumed to be known exactly" on the accuracy of the state estimation algorithm and eventually the optimal placement of measurement devices, in case line parameters are uncertain.

ACKNOWLEDGEMENT

This research project is financially supported by the Swiss Innovation Agency Innosuisse and is part of the Swiss competence Center for Energy Research SCCER FURIES. The authors would also like to acknowledge the technical and financial support of "Services Industriels de Genève - SIG", utility of Geneva, Switzerland, in the frame of the POEM project No. 88824.